\documentclass{article}
\usepackage{}

\usepackage{amssymb,amsmath}
\usepackage{latexsym}
\usepackage[dvips]{graphicx}
\usepackage{graphicx,color}
\usepackage{hyperref}

\newtheorem{theorem}[subsection]{Theorem}
\newtheorem{lemma}[subsection]{Lemma}

\newcommand{\qed}{\hfill\ensuremath{\square}}

\begin{document}
\title{On cross-$2$-intersecting families \thanks{E-mail addresses: 23787267@qq.com (Yanhong Chen), anshuili@usx.edu.cn (Anshui Li),
wu@hunnu.edu.cn (Biao Wu), huajunzhang@usx.edu.cn (Huajun Zhang, Corresponding author)}}
\date{}

\author{
Yanhong Chen $^{a}$,
Anshui Li $^{b}$,
Biao Wu $^{c}$
Huajun Zhang $^{a}$
 \\[2ex]
{\small $^{a}$ Department of Mathematics, Shaoxing University} \\
{\small Shaoxing, Zhejiang, 321004, P.R. China }\\
{\small $^{b}$ Department of Statistics, Shaoxing University} \\
{\small Shaoxing, Zhejiang, 321004, P.R. China }\\
{\small $^{c}$ MOE-LCSM,  School of Mathematics and Statistics, Hunan Normal University} \\
{\small Changsha, Hunan, 410081, P.R. China } }

\maketitle

\begin{abstract}
 Two families $\mathcal A\subseteq\binom{[n]}{k}$ and $\mathcal B\subseteq\binom{[n]}{\ell}$ are called cross-$t$-intersecting if $|A\cap B|\geq t$ for all  $A\in\mathcal A$, $B\in\mathcal B$.
Let $n$, $k$ and $\ell$ be positive integers such that $n\geq 3.38\ell$ and $\ell\geq k\geq 2$. In this paper, we will determine the upper bound of $|\mathcal A||\mathcal B|$ for cross-$2$-intersecting families $\mathcal A\subseteq\binom{[n]}{k}$ and $\mathcal B\subseteq\binom{[n]}{\ell}$. The structures of the extremal families attaining the upper bound are also characterized. The similar result obtained by Tokushige can be considered as a special case of ours when $k=\ell$, but under a more strong condition $n>3.42k$. Moreover, combined with the results obtained in this paper, the complicated extremal structures attaining the upper bound for nontrivial cases can be relatively easy to reach with similar techniques.
\end{abstract}



\section{Introduction}

Let $n$, $k$  and $t$ be positive integers with $t\leq k\leq n$.  We write $[n]:=\{1,\ldots,n\}$ and denote its power set by $2^{[n]}$. A nonempty subset $\mathcal{F}$ of $2^{[n]}$ is called {\it $k$-uniform} if its elements all have size $k$.  We use $\binom{[n]}{k}$
to denote the set of $\{A\subseteq [n]:|A|=k\}$. We call $\mathcal{F}$ a {\it $t$-intersecting family} if $|A\cap B|\ge t$ for all $A,B\in \mathcal{F}$. A $1$-intersecting family is simply called an {\it intersecting} family. The following Erd\H{o}s--Ko--Rado theorem is one of the central results in extremal combinatorics.

\begin{theorem}[Erd\H{o}s--Ko--Rado \cite{EKR1961}]\label{EKR1961}
Let $t,k,n$ be positive integers such that $t\leq k\leq n$. Suppose that  $\mathcal F$ is a $t$-intersecting family of $\binom{[n]}{k}$. Then for $n\geq n_0(k,t)$, \[|\mathcal{F}|\leq \binom{n-t}{k-t}. \]
\end{theorem}

Let $k,t$ be fixed positive integers  such that $t\le k$, and let $N_0(k,t)$ be the smallest possible value of $n_0(k,t)$ in Theorem \ref{EKR1961}. In the case $t=1$, it is shown in \cite{EKR1961} that $N_0(k,1)=2k$. For a general $t$ we have $N_0(k,t)=(t+1)(k-t+1)$, which is established in \cite{Fra1976} for $t\ge 15$ and in \cite{Wil1984} for all $t$.

Let $\mathcal{A}$ and $\mathcal{B}$ be two families of subsets of $[n]$. They are called {\it cross-$t$-intersecting} if $|A\cap B|\ge t$ for all $A\in \mathcal{A}$ and $B\in \mathcal{B}$.
In the case $t=1$, we simply say that $\mathcal{A}$ and $\mathcal{B}$ are {\it cross-intersecting}.
The cross-$t$-intersecting property is a natural extension of the $t$-intersecting property, since the two properties coincide in the case $\mathcal{A}=\mathcal{B}$. Let $\mathcal A$ and $\mathcal B$ be cross-$t$-intersecting families of $\binom{[n]}{k}$. We say that they are \emph{maximal} (or \emph{maximum}) if there are no cross-$t$-intersecting families $\mathcal A_1$ and $\mathcal B_1$ of $\binom{[n]}{k}$ such that $\mathcal A\subseteq \mathcal A_1$,  $\mathcal B\subseteq\mathcal B_1$ and  $|\mathcal A||\mathcal B|<|\mathcal A_1||\mathcal B_1|$ (or
{$|\mathcal A||\mathcal B|<|\mathcal A_1||\mathcal B_1|$}).

Let $\mathcal A\subseteq\binom{[n]}{k}$ and $\mathcal B\subseteq\binom{[n]}{\ell}$ be two cross-intersecting families. In 1986, Pyber \cite{Pyber} determined the upper bound under the assumption that $n\ge 2k+\ell-2$ and $k\geq \ell$ or $k=\ell$ and $n\ge 2k$, and then  Matsumoto and Tokushige \cite{MT1989} proved that the assumption can be improved to be $n\geq\max\{2k,2\ell\}$.

\begin{theorem}[ \cite{MT1989,Pyber}] \label{theoMT}
If $\mathcal{A}\subseteq {\binom{[n]}{k}}$ and $\mathcal{B}\subseteq {\binom{[n]}{\ell}}$ are cross-intersecting with $n\ge \max\{2k,2\ell\}$, then
$$|\mathcal{A}||\mathcal{B}|\leq {\binom{n-1}{k-1}}{\binom{n-1}{\ell-1}}.$$ Moreover, the equality holds if and only if $\mathcal{A}=\{A\in {\binom{[n]}{k}}: x\in A\}$ and $\mathcal{B}=\{B\in \binom{[n]}{\ell}: x\in B\}$ for some fixed element $x\in[n]$ unless $n=2k=2\ell$.
\end{theorem}

In 2013, Tokushige \cite{Tokushige2013} proved the following analogous result for cross-$t$-intersecting families.

\begin{theorem} [Tokushige \cite{Tokushige2013}] \label{theoTokushige}
Suppose that $n,k,t$ are natural numbers such that $t\le k\le n$ and $\frac{k}{n}<1-\frac{1}{\sqrt[t]{2}}$. If two $k$-uniform families $\mathcal{A}$, $\mathcal{B}$ are cross-$t$-intersecting, then $ |\mathcal{A}|\cdot|\mathcal{B}|\leq \binom{n-t}{ k-t}^2.$ Moreover, the equality holds if and only if $\mathcal{A}=\mathcal{B}=\{F\in {\binom{[n]}{k}}: T\subset F\}$ for some $t$-subset $T$ of $[n]$.
\end{theorem}

In the same paper \cite{Tokushige2013}, Tokushige conjectures that the lower bound on $n$ in Theorem \ref{theoTokushige} can be improved to $n \ge (t+1)(k-t+1)$. In 2014, Frankl, Lee,  Siggers and Tokushige \cite{FLST} verified the conjecture under the stronger assumptions that $t\ge 14$ and $n \ge (t+1)k$. In the same year, Borg \cite{Borg2014} independently confirmed the conjecture for large $n$. Recently, Zhang and Wu \cite{ZW} confirmed Tokushige's conjecture for all $t\ge 3$.


 In the present paper, we will consider the case of $t=2$ and obtain the following result.

\begin{theorem} \label{mainthm}If $\mathcal{A}\subseteq {\binom{[n]}{k}}$ and $\mathcal{B}\subseteq {\binom{[n]}{\ell}}$ are cross-$2$-intersecting with $n\ge 3.38\max\{k,\ell\}$, then $$|\mathcal A||\mathcal B|\leq \binom{n-2}{k-2}\binom{n-2}{\ell-2}.$$
Moreover, the equality holds if and only if $\mathcal A=\{A\in\binom{[n]}{\ell}: T\subset A\}$ and $\mathcal B=\{B\in\binom{[n]}{k}: T\subset B\}$ for some $2$-subset $T$ of $[n]$.
\end{theorem}

The methods used in this paper are the shift operator method and the generating set method.
Their definitions and basic properties will be given in Section~\ref{sec:pre}. In Section~\ref{sec:inq}, we will prove two inequalities concerning the proof of Theorem~\ref{mainthm}. The whole proof of Theorem~\ref{mainthm} will be presented in Section~\ref{sec:main}.
The complicated extremal structures attaining the upper bound for nontrivial cases are also shown as Theorem~\ref{thm:non} without proof in Section~\ref{sec:nontrivial}.
\section{Preliminaries}\label{sec:pre}
In the study of Erd\H{o}s-Ko-Rado theorem, the shift operator is a powerful method which can be traced back to \cite{EKR1961}.  For $\mathcal A\subset\binom{[n]}{k}$ and $1\leq i<j\leq n$,  define
\[
  s_{ij}(A)=\begin{cases}(A\setminus\{j\})\cup\{i\}, & \hbox{if $j\in A,\ i\not\in A, (A\setminus\{j\})\cup\{i\}\not\in \mathcal{A}$;}\\A, & \hbox{otherwise,}\end{cases}
\]
and set $s_{ij}(\mathcal A)=\{s_{ij}(A):A\in\mathcal A\}$ correspondingly. It is well known that the shift operator preserving the cross-$t$-intersecting property \cite{Fra1987}.
A family $\mathcal A$ is called \emph{left-compressed} if $s_{ij}(\mathcal A)=\mathcal A$ for all $1\leq i<j\leq n$.

For $\mathcal E\subseteq 2^{[n]}$, we say $\mathcal E$ an \emph{upset} if $F\in \mathcal E$ whenever $E\subseteq F$ for some $E\in \mathcal E$. Let $\mathcal A$ be an intersecting family of $\binom{[n]}{k}$. An upset $\mathcal E$ is called a \emph{generating set}, which was first introduced in \cite{AK1996}, of $\mathcal A$ if $\mathcal E\cap \binom{[n]}{k}=\mathcal A$. Clearly, $\mathcal A$ is a generating set of itself. Then all the generating sets of $\mathcal A$ form a nonempty set denoted by $G(\mathcal A)$. For $g(\mathcal A)\in G(\mathcal A)$, let $g_*(\mathcal A)$ be the set of all minimum elements under inclusion. Set $G_*(\mathcal A)=\{g_*(\mathcal A): g(\mathcal A)\in G(\mathcal A)\}$.

For $E\in 2^{[n]}$ and $\mathcal E\subseteq 2^{[n]}$, set $s^+(E)=\max\{i: i\in E\}$ and $s^+(\mathcal E)=\max\{s^+(E): E\in\mathcal E\}$. For $g(\mathcal A)\in G_*(\mathcal A)$, set $s=s^+(g(\mathcal A))$, let $g^*(\mathcal A)=\{E\in g(\mathcal A): s\in E\}$ and $g^*_i(\mathcal A)=\{E\in g^*(\mathcal A): |E|=i\}$ for $t\leq i\leq s$. For $\mathcal E\subseteq g^*(\mathcal A)$, set $\mathcal E'=\{E\setminus\{s\}: E\in\mathcal E\}$.  From \cite{AK}, we know that the generating sets of $\mathcal A$ have the following properties.

\begin{lemma}\label{gtf} Let $n$, $k$ and $t$ be positive integers such that $n>2k-t$, and let  $\mathcal A$ be a left-compressed $t$-intersecting family of $\binom{[n]}{k}$. For $g(\mathcal A)\in G_*(\mathcal A)$, set  $s=s^+(g(\mathcal A))$, then

  {\rm (i)} $g(\mathcal A)$ is $t$-intersecting.

  {\rm (ii)} For $1\leq i<j\leq s$  and $E\in g(\mathcal A)$, $F\subseteq s_{ij}(E)$ for some $F\in g(\mathcal A)$.

  {\rm(iii)} $\mathcal A$ is a disjoint union \[\mathcal A=\bigcup_{E\in g(\mathcal A)}\mathcal D(E), \mbox{ where }
  \mathcal D(E)=\{B\in\binom{[n]}{k}: B\cap [s^+(E)]=E\}.\]

  {\rm(iv)} If $\mathcal A$ is maximal, then for any $E_1, E_2\in g^*(\mathcal A)$ with $|E_1\cap E_2|=t$, necessarily $|E_1|+|E_2|=s+t$ and $E_1\cup E_2=[s]$. That is, if $g^*_i(\mathcal A)\neq \emptyset$, then $g^*_{s+t-i}(\mathcal A)\neq \emptyset$ and for any $E_1\in g^*_i(\mathcal A)$, there exists $E_2\in g^*_{s+t-i}(\mathcal A)$ with $|E_1\cap E_2|=t$ and $E_1\cup E_2=[s]$.

{\rm(v)} If  $g^*_i(\mathcal A)\neq\emptyset$, $\mathcal A_1=\mathcal A\cup\mathcal D(g^{*}_i(\mathcal A)')\backslash \mathcal D(g^{*}_{s+t-i}(\mathcal A))$ is also $t$-intersecting family of $\binom{[n]}{k}$ with
$$|\mathcal A_1|=|\mathcal A|+|g^{*}_i(\mathcal A)|\binom{n-s}{k-i+1}-|g^{*}_{s+i-t}(\mathcal A)|\binom{n-s}{k+i-s-t}.$$
 \end{lemma}
Suppose that  $\mathcal A\subseteq \binom{[n]}{k}$ and $\mathcal B\subseteq\binom{[n]}{\ell}$. For $n>k+\ell-t$, it is easy to find that $\mathcal A$  and $\mathcal B$ are cross-$t$-intersecting if and only if $g(\mathcal A)$ and $g(\mathcal B)$ are so for all $g(\mathcal A)\in G(\mathcal A)$ and $g(\mathcal B)\in G(\mathcal B)$.  Similarly, on the generating sets of $\mathcal A$ and $\mathcal B$, we also have the following properties.

\begin{lemma}\label{lemmag1}
 Let $n$, $k$, $\ell$ and $t$ be positive integers with $n> k+\ell-t$, and let $\mathcal A\subseteq\binom{[n]}{k}$ and $\mathcal B\subseteq\binom{[n]}{\ell}$ be maximal left-compressed cross-$t$-intersecting families. For  $g(\mathcal A)\in G_*(\mathcal A)$ and $g(\mathcal B)\in G_*(\mathcal B)$. Set $s=\max\{s^+(g(\mathcal A)), s^+(g(\mathcal B))\}$. Then,
\begin{itemize}
\item[\rm(i)]For $1\leq i<j\leq s$  and $E\in g(\mathcal A)$ (or $g(\mathcal B)$), we have either $s_{ij}(E)\in g(\mathcal F)$ or $F\subset s_{ij}(E)$ for some $F\in g(\mathcal A)$ (or $g(\mathcal B)$).
    \item [\rm(ii)] For $t\leq i\leq k$,  $g_i^*(\mathcal A)\neq\emptyset$ if and only if $g_{s+t-i}^*(\mathcal B)\neq\emptyset$. Furthermore, for each $E\in g_i^*(\mathcal A)\neq\emptyset$, there exists $F\in g_{s+t-i}^*(\mathcal B)$ such that $|E\cap F|=t$, $E\cup F=[s]$.
        \item[\rm(iii)]If $g_i^*(\mathcal A)\neq\emptyset$, $\mathcal A_1=\mathcal A\cup \mathcal{D}(g_i^*(\mathcal A)')$ and $\mathcal B_1=\mathcal B\setminus \mathcal D(g_{s+t-i}^*(\mathcal B))$ are also cross-$t$-intersecting families of $\binom{[n]}{k}$ with
            \[|\mathcal A_1|=|\mathcal A|+|g_i^*(\mathcal A)|\binom{n-s}{k-i+1} \ and \ |\mathcal B_1|=|\mathcal B|-|g_{s+t-i}^*(\mathcal B)|\binom{n-s}{\ell+i-s-t}. \]

 \end{itemize}
\end{lemma}
\textbf{Proof.} (i) and (iii) are obvious. So we only need to prove (ii). Suppose $g_i^*(\mathcal A)\neq\emptyset$ for some $t\leq i\leq k$. For any $E\in g_i^*(\mathcal A)$, set $E'=E\setminus\{s\}$, it is clear that  $\mathcal D(E')\not\subseteq \mathcal A$. If $|E'\cap F|\geq t$ for all $F\in g(\mathcal B)$, then $\mathcal A_1=\mathcal A\cup \mathcal D(E')$ and $\mathcal B$ are also cross-$t$-intersecting with $\mathcal A\subsetneq\mathcal A_1$, contradicting the maximality of $\mathcal A$ and $\mathcal B$. Therefore, there exists $F_1\in g(\mathcal B)$ such that $|E'\cap F_1|=t-1$. On the other hand, $|E\cap F_1|\geq t$, hence $|E\cap F_1|=t$ and $s\in E\cap F_1$. If $r\not\in E\cup F_1$ for some $r\in [s]$, then by (i) there exists $F_2\in g(\mathcal B)$ such that $F_2\subseteq s_{rs}(F_1)$, but $s\not\in F_2\cap E\subsetneq F_1\cap E$. Consequently, $|E\cap F_2|<t$, contradicting that $g(\mathcal A)$ and $g(\mathcal B)$ are cross-$t$-intersecting. Hence, $E\cup F_1=[s]$ and $|E\cap F_1|=t$, i.e, $|F_1|=s+t-i$ and so $g_{s+t-i}^*(\mathcal B)\neq\emptyset$. This completes the proof.
\qed



\section{Some basic inequalities concerning the proofs}\label{sec:inq}

For $\max\{s+2-k,3\}\leq i\leq \min\{s-1, k\}$,
set
\begin{align*}f(n,k,s,i)&=
\frac{\binom{s-1}{i-1}\binom{n-s}{k-i}+\binom{s-1}{i}\binom{n-s+1}{k-i}}
{\binom{s-1}{i-1}\binom{n-s+1}{k-i+1}+\binom{s-1}{i}\binom{n-s+1}{k-i}}\\
&=\frac{(k-i+1)(s(n-s+1)-i(k-i))}{(n-s+1)(i(n-k-s+i+1)+(s-i)(k-i+1))}\\
&=\frac{(k-i+1)(s(n-s+1)-i(k-i))}{(n-s+1)(i(n+2i-2k-2s)+s(k+1))}.
\end{align*}
\begin{lemma}\label{fksi}
Let $n,k,s$ and $i$ be positive integers with $n\geq 3.38k$ and
$$\max\{s+2-k,3\}\leq i\leq \min\{s-1, k\}.$$ If $s\geq k+1$, then
\[f(n,k,s,i)>\frac{\binom{n-s}{k-i}}{\binom{n-s}{k-i+1}}.\]
\end{lemma}
Proof.  By definition, it is easy to find that \begin{eqnarray*}f(n,k,s,i)&=&\frac{\binom{s-1}{i-1}\binom{n-s}{k-i}+\binom{s-1}{i}\binom{n-s+1}{k-i}}
{\binom{s-1}{i-1}\binom{n-s}{k-i+1}+\binom{s-1}{i}\binom{n-s+1}{k-i}+\binom{s-1}{i-1}\binom{n-s}{k-i}}>\frac{\binom{n-s}{k-i}}{\binom{n-s}{k-i+1}}=\frac{k-i+1}{n+i-k-s}
\end{eqnarray*}
if and only if
 $$\binom{s-1}{i}\binom{n-s+1}{k-i}(1-\frac{k-i+1}{n+i-k-s})>\binom{s-1}{i-1}\binom{n-s}{k-i}\frac{k-i+1}{n+i-k-s},$$
which is equivalent to
\begin{eqnarray}\label{eq1}(n-s+1)(s-i)(n-2k-s+2i-1)-i(n-s-k+i+1)(k-i+1)>0.\end{eqnarray}
The assumption $\max\{s+2-k,3\}\leq i\leq \min\{s-1, k\}$ implies $s+2\leq k+i$ and $i\leq k$, it follows that   $n-s+1\geq n-s-k+i+1$ and $n-2k-s+2i-1\geq n-3k+i+1>i$, and the two equalities both hold if and only if $k=i$, $s=2k-1$ and $n=3.38k$. Hence,  the inequality (\ref{eq1}) holds whenever $s\geq k+1$.

\begin{lemma}\label{lem:main}
Let $n$, $k$, $\ell$ and $s$ be four positive integers with $n\geq 3.38\ell$, $s\leq k+\ell-2$ and $\max\{s+2-\ell,3\}\leq i\leq \min\{k,s-1\}$. Then
\begin{equation}\label{inq:main}
f(n,k,s,i)f(n,\ell,s,s+2-i)>\frac{\binom{n-s}{k-i}\binom{n-s}{\ell+i-s-2}}{\binom{n-s}{k-i+1}\binom{n-s}{\ell+i-s-1}}
\end{equation}unless $(s,i)=(4,3)$, $(5,3)$, $(5,4)$ and $(6,4)$.\end{lemma}
\textbf{Proof.} By Lemma~\ref{fksi}, the result holds if $s> \max\{k,\ell\}$. In the following, we always assume that $\max\{k,\ell\}\geq s$.
By definition, to complete the proof, we only need to prove that $T(n,s,k,\ell,i) >1$
unless $(s,i)=(4,3)$, $(5,3)$, $(5,4)$ and $(6,4)$,
where
\begin{small}
\begin{eqnarray*}
T(n,s,k,\ell,i) =
\frac{(s(n-s+1)-i(k-i))(s(n-s+1)-(s-i+2)(\ell+i-s-2))(n-s-k+i)(n-\ell-i+2)}
{(i(n-2k-s+2i)+s(k-i+1))((s-i+2)(n-2\ell-2i+s+4)+s(\ell+i-s-1))(n-s+1)^2}.
\end{eqnarray*}
\end{small} We distinguish three cases to complete the proof.

{{Case 1: $5\leq i\leq s-3$}}. Clearly, $5\leq s+2-i\leq s-3$.
Set
 $$g(n,s,k,\ell,i)=(n-s-k+i)(s(n-s+1)-(s-i+2)(\ell+i-s-2))-(n-s+1)(i(n-2k-s+2i)+s(k-i+1))$$
 and
 $$h(n,s,k,\ell,i)=(n-\ell-i+2)(s(n-s+1)-i(k-i))-(n-s+1)((s-i+2)(n-2\ell-2i+s+4)+s(\ell+i-s-1)).$$
 To complete the proof, it suffices to verify that
 \[g(n,s,k,\ell,i)>0 \:\:\text{and}\:\:  h(n,s,\ell,k,i)>0. \]
 By definition,
 \begin{eqnarray*}
 g'_k(n,s,k,\ell,i)&=&-(s(n-s+1)-(s+2-i)(\ell+i-s-2))-(n-s+1)(s-2i)\\&=&n(-2s+2i)+\ell(s-i+2)+s^2-6s-i^2+6i-4\\&\leq& \ell(-5.76s+5.76i+2)+s^2-6s-i^2+6i-4\\&\leq& s(-4.76s+5.76i-4)-i^2+6i-4\\&\leq& -5.52i-31.04<0
 \end{eqnarray*}
and $$g'_{\ell}(n,s,k,\ell,i)=-(n-s-k+i)(s+2-i)<0.$$
It follows that $g(n,s,k,\ell,i)$ is decreasing function in both $k$ and $\ell$.  Noting that $h(n,s,k,\ell,i)=g(n,s,\ell,k,s+2-i)$, so we only need to prove that $g(n,s,\ell,\ell,i)>0$ for $5\leq i\leq s-3$ under the assumption that $n \geq 3.38\ell$.
As $s\geq i+2$,  $g''_n(n,s,\ell,\ell,i)=2s-2i\geq0$, it follows that $g'_n(n,s,\ell,\ell,i)$ is an increasing function in $n$.
Hence \begin{eqnarray*}
g'_n(n,s,\ell,\ell,i)&\geq&g'_n(3.38\ell,s,\ell,\ell,i)\\&=&\ell(5.76s-5.76i-2)-i^2+2is-5i-s^2+4s+4\\&\geq&s(4.76s-3.76i+2)-i^2-5i+4\\&\geq&8.52i+27.04>0
\end{eqnarray*}
subject to the assumption that  $\ell\geq s\geq i+2$, $i\geq5$ and $n\geq 3.38\ell$. Consequently, we conclude that  $g(n,s,\ell,\ell,i)$ is an increasing function with respect to $n$.
Moreover,   $g''_s(3.38\ell,s,\ell,\ell,i)=-2.76\ell-8<0$,  it follows that $g(3.38\ell,s,\ell,\ell,i)$ is a
concave function in $s$ in $[i+3,\ell]$. Hence,$$g(3.38\ell,s,\ell,\ell,i) \geq \min\{g(3.38\ell,i+3,\ell,\ell,i),g(3.38\ell,\ell,\ell,\ell,i)\}$$ for $\ell-3\geq i\geq 5$.  However, for $i\geq 5$ and $\ell\geq s\geq i+3$, $g'_\ell(3.38\ell,\ell,\ell,\ell,i)=2.712\ell^2+1.52\ell+4.952i\ell-4.38i^2-1.9i+4.52\geq 3.284i^2+20.372i+18.408>0$ and $g'_l(3.38\ell,i+3,\ell,\ell,i)=9.28\ell-i+7>0$,  so we have
\begin{eqnarray*}
g(3.38\ell,\ell,\ell,\ell,i)&\geq&g(3.38i+10.14,i+3,i+3,i+3,i)=2.7132i^2+54.0792i+44.8188>0
\end{eqnarray*}
and \begin{eqnarray*}g(3.38\ell,i+3,\ell,\ell,i)&\geq&g(3.38i+10.14,i+3,i+3,i+3,i)=2.7132i^2+54.0792i+44.8188>0.\end{eqnarray*}
Therefore, $T(n,s,k,\ell,i)>1$ for $5\leq i \leq s-3$.

Case 2: $i=4$ and $s\geq 7$. Then
\begin{align*}
 T(n,s,k,\ell,4)=
{\tiny \frac{(s(n-s+1)-4(k-4))(s(n-s+1)-(s-2)(\ell-s+2))(n-s-k+4)(n-\ell-2)}
{(4(n-2k-s+8)+s(k-3))((s-2)(n-2\ell+s-4)+s(\ell-s+3))(n-s+1)^2}.}
\end{align*}

We first consider the case  $\ell \geq k$.  Note that
\begin{eqnarray*}
&&((s-2)(n-2\ell+s-4)+s(\ell-s+3))-0.38\ell-(4(n-2k-s+8)+s(k-3))\\
&\geq&((s-2)(3.38\ell-2\ell+s-4)+s(\ell-s+3))-0.38\ell-(4(3.38\ell-2k-s+8)+s(k-3))\\
&\geq&1.38ks-8.66k+4s-24>0.
\end{eqnarray*}
Hence
\begin{small}
\begin{eqnarray*}
T(n,s,k,\ell,4)
\geq \frac{(s(n-s+1)-4(k-4))(s(n-s+1)-(s-2)(\ell-s+2))(n-s-k+4)(n-\ell-2)}
{(4(n-2k-s+8)+s(k-3)+0.38\ell)((s-2)(n-2\ell+s-4)+s(\ell-s+3)-0.38\ell)(n-s+1)^2}.
\end{eqnarray*}
\end{small}
Set
\begin{small}
\begin{eqnarray*}
p(n,s,k,\ell)
=(s(n-s+1)-(s-2)(\ell-s+2))(n-s-k+4)-(4(n-2k-s+8)+s(k-3)+0.38\ell)(n-s+1)
\end{eqnarray*}
\end{small}
and
\begin{small}
\begin{eqnarray*}
q(n,s,k,\ell)=(s(n-s+1)-4(k-4))(n-\ell-2)-((s-2)(n-2\ell+s-4)+s(\ell-s+3)-0.38\ell)(n-s+1).
\end{eqnarray*}
\end{small}
As
 \begin{eqnarray*}
p'_n(n,s,k,\ell)&=&(2s-8)n-2ks+8k-s\ell+1.62\ell-s^2+12s-32\\
&\geq&5.76s\ell-25.42\ell-s^2+12s-2ks+8k-32\\
&\geq&3.76s\ell-17.42\ell-s^2+12s-32>0
\end{eqnarray*}
and
\begin{eqnarray*}
q'_n(n,s,k,\ell)=4n-4k-3.62\ell-s+10>0,
\end{eqnarray*}
it follows that   $p(n,s,k,\ell)\geq p(3.38k,s,k,\ell)$ and $q(n,s,k,\ell)\geq q(3.38k,s,k,\ell)$. Moreover,  $p''_s(3.38k,s,k,\ell)=-4.76\ell+2k-8<0$ and $q''_s(3.38k,s,k,\ell)=-2<0$.
Similarly, to complete the proof, it remains to verify that
\[p(3.38\ell,\ell,k,\ell)\geq 0, p(3.38\ell,7,k,\ell)\geq 0, q(3.38\ell,\ell,k,\ell)\geq 0 \:\: \text{and} \:\: q(3.38\ell,7,k,\ell)\geq 0.\]
 Noting
$p'_{k}(n,s,k,\ell)=-s(n-s+1)+(s-2)(\ell-s+2)+8<0$ and $q'_{k}(n,s,k,\ell)=-4(n-\ell+2)<0$,
 by the assumptions that $\ell\geq s\geq 7$ and $\ell\geq k$, it follows that
\begin{eqnarray*}
p(3.38\ell,\ell,k,\ell)&\geq &p(3.38\ell,\ell,\ell,\ell)
=0.9044\ell^3+9.758\ell^2-73.54\ell-16>0,
\end{eqnarray*}
\begin{eqnarray*}
p(3.38\ell,7,k,\ell)\geq p(3.38\ell,7,\ell,\ell)=0.8088\ell^2+38.42\ell-51>0
\end{eqnarray*}
and
\begin{eqnarray*}
q(3.38\ell,\ell,k,\ell)\geq q(3.38\ell,\ell,\ell,\ell)=0.3332\ell^2+31.18\ell-40>0,
\end{eqnarray*}
\begin{eqnarray*}
q(3.38\ell,7,k,\ell)\geq q(3.38\ell,7,\ell,\ell)=1.0932\ell^2+23.86\ell-26>0.
\end{eqnarray*}
Hence $T(n,s,k,\ell,4)>1$ if $\ell \geq k$.

Now we  consider the case $k\geq \ell$. Since
$
(n-\ell-2)+0.1k-(n-s-k+4)>0,
$
we have
\begin{eqnarray*}
&&T(n,s,k,\ell,4)\\
&\geq& \frac{(s(n-s+1)-4(k-4))(s(n-s+1)-(s-2)(\ell-s+2))(n-s-k+4-0.1k)(n-\ell-2+0.1k)}
{(4(n-2k-s+8)+s(k-3))((s-2)(n-2\ell+s-4)+s(\ell-s+3))(n-s+1)^2}.
\end{eqnarray*}
Set
\begin{eqnarray*}
&&s(n,s,k,\ell)\\
&=&(s(n-s+1)-(s-2)(\ell-s+2))(n-s-k+4-0.1k)-(4(n-2k-s+8)+s(k-3))(n-s+1)
\end{eqnarray*}
and
\begin{eqnarray*}
&&t(n,s,k,\ell)\\
&=&(s(n-s+1)-4(k-4))(n-\ell-2+0.1k)-((s-2)(n-2\ell+s-4)+s(\ell-s+3))(n-s+1).
\end{eqnarray*}
By definition and calculation, we can also obtain that $s'_n(n,s,k,\ell)>0$, $t'_n(n,s,k,\ell)>0$, $s'_{\ell}(n,s,k,\ell)<0$, $t'_{\ell}(n,s,k,\ell)<0$, $s''_s(n,s,k,\ell)<0$ and $t''_s(n,s,k,\ell)<0$.

Similarly, to prove that $T(n,s,k,\ell,4)>0$, it remains  to verify that
 $s(3.38k,k,k,\ell)\geq 0$ , $s(3.38k,7,k,\ell)\geq 0$ and $t(3.38k,k,k,\ell)\geq 0$ , $t(3.38k,7,k,\ell)\geq 0$. However, by the assumption that $k\geq s\geq 7$ and $k\geq \ell$, it follows that
\begin{eqnarray*}
s(3.38k,k,k,\ell)&\geq&s(3.38k,k,k,k)=
0.8806k^3+10.6724k^2-73.2k-16>0,
\end{eqnarray*}
\begin{eqnarray*}
s(3.38k,7,k,\ell)&\geq&s(3.38k,7,k,k)=1.9066k^2+36.31k-51>0
\end{eqnarray*}
and
\begin{eqnarray*}
t(3.38k,k,k,\ell)&\geq&t(3.38k,k,k,k)=
0.0238k^3-0.6012k^2+30.96k-40>0,
\end{eqnarray*}
\begin{eqnarray*}
t(3.38k,7,k,\ell)&\geq&t(3.38k,7,k,k)=0.0054k^2+25.88k-26>0.
\end{eqnarray*}
Hence $T(n,s,k,\ell,4)>1$ for $k \geq \ell$.

Therefore, $T(n,s,k,\ell,i)>1$ subject to  $i=4$ and $s \geq 7$.

Case 3: $i=3$.

Suppose $s\geq 7$. Then
\begin{eqnarray*}
T(n,s,k,\ell,3)=\frac{(s(n-s+1)-(s-1)(\ell-s+1))(s(n-s+1)-3(k-3))(n-\ell-1)(n-s-k+3)}
{((s-1)(n-2\ell+s-2)+s(\ell-s+2))(3(n-2k-s+6)+s(k-2))(n-s+1)^2}
\end{eqnarray*}
If $\ell \geq k$, since $(n-s+1-0.27\ell)(n-s+1+0.31\ell)-(n-s+1)(n-s+1)\geq0.0515\ell^2-0.04s\ell+40\ell>0$,
 we have
 {\small
 \begin{eqnarray*}
&&~~T(n,s,k,\ell,3)\\
&>&\frac{(s(n-s+1)-(s-1)(\ell-s+1))(s(n-s+1)-3(k-3))(n-\ell-1)(n-s-k+3)}
{((s-1)(n-2\ell+s-2)+s(\ell-s+2))(3(n-2k-s+6)+s(k-2))(n-s+1-0.27\ell)(n-s+1+0.31\ell)}.
\end{eqnarray*}}
Set
\begin{eqnarray*}
&&p(n,s,k,\ell)\\
&=&(s(n-s+1)-3(k-3))(n-\ell-1)-((s-1)(n-2\ell+s-2)+s(\ell-s+2))(n-s+1-0.27\ell)
\end{eqnarray*}
and
\begin{eqnarray*}
&&q(n,s,k,\ell)\\
&=&(s(n-s+1)-(s-1)(\ell-s+1))(n-s-k+3)-(3(n-2k-s+6)+s(k-2))(n-s+1+0.31\ell).
\end{eqnarray*}
By the assumption $\ell\geq s\geq 7$,  it follows that
\begin{eqnarray*}
 p(3.38\ell,s,k,\ell)&=&\ell(0.6426s\ell+4.2918\ell-7.14k-1.65s)+3k+16.58\ell+20s-11\\
&\geq&\ell(8.79\ell-7.14k-1.65s)+3k+16.58\ell+20s-11\\
&\geq&\ell(8.79\ell-7.14\ell-1.65\ell)+3k+16.58\ell+20s-11\\
&=&3k+16.58\ell+20s-11>0,
\end{eqnarray*}
\begin{eqnarray*}
q(3.38\ell,\ell,k,\ell)&=&k(-5.07\ell^2+15.14\ell+5)+5.664\ell^3-6.687\ell^2-51.18\ell-15\\
&\geq&0.5944\ell^3+8.4534\ell^2-46.18\ell-15>0
\end{eqnarray*}and
\begin{eqnarray*}
q(3.38\ell,7,k,\ell)=\ell(22.2742\ell-21.35k)+12k+32.65\ell-78>0.
\end{eqnarray*}
 Noting that $p'_n(n,s,k,\ell)=2n-3k-2.27\ell+0.27s\ell-s+8>0$, $q'_n(n,s,k,\ell)=\ell(22.2742\ell-21.35k)+12k+32.65\ell-78>0$ and $q''_s(3.38\ell,s,k,\ell)=2k-4.76\ell-8<0$, we conclude that $p(n,s,k,\ell)>0$ and $q(n,s,k,\ell)\geq \max\{q(3.38\ell,7,k,\ell), q(3.38\ell,\ell,k,\ell)\}>0$.

If $\ell<k$,
 since $(n-\ell-1)+0.12k-(n-s-k+3)>0,$ we then have
 {
\begin{eqnarray*}
&&T(n,s,k,\ell,3)\\
&\geq&
\frac{(s(n-s+1)-(s-1)(\ell-s+1))(s(n-s+1)-3(k-3))(n-\ell-1+0.12k)(n-s-k+3-0.12k)}
{((s-1)(n-2\ell+s-2)+s(\ell-s+2))(3(n-2k-s+6)+s(k-2))(n-s+1)^2}.
\end{eqnarray*}
}
Set
\begin{eqnarray*}
&&p(n,s,k,\ell)\\
&=&(s(n-s+1)-3(k-3))(n-\ell-1+0.12k)-((s-1)(n-2\ell+s-2)+s(\ell-s+2))(n-s+1)
\end{eqnarray*}
and
\begin{eqnarray*}
&&q(n,s,k,\ell)\\
&=&(s(n-s+1)-(s-1)(\ell-s+1))(n-s-k+3-0.12k)-(3(n-2k-s+6)+s(k-2))(n-s+1).
\end{eqnarray*}
Noting that
$$p'_n(n,s,k,\ell)=2n-3k-2\ell-s+0.12ks+8>0,$$
and
$$q'_n(n,s,k,\ell)=2sn-6n-2.12ks+6k-s\ell+\ell+10s-20\geq k(4.64s-14.28)-s\ell+\ell+10s-20>0.$$

 It is easy to see that $p(n,s,k,\ell)\geq p(3.38k,s,k,\ell)$ and $q(n,s,k,\ell)\geq q(3.38k,s,k,\ell)$. Moreover, $p''_s(3.38k,s,k,\ell)=-0.24k<0$ and $q''_s(3.38k,s,k,\ell)=-4.76k+2\ell-8<0$. To complete the proof, it suffices to verify $p(3.38k,k,k,\ell)\geq 0$, $p(3.38k,7,k,\ell)\geq 0$ and $q(3.38k,k,k,\ell)\geq 0$, $q(3.38k,7,k,\ell)\geq 0$. However, by the assumption that $k\geq \ell$, we have
\begin{eqnarray*}
p(3.38k,k,k,\ell)=0.2856k^3-2.3356k^2-1.76k\ell+33.12k-11\ell-11\geq0.2856k^3-4.0956k^2+22.12k-11,
\end{eqnarray*}
\begin{eqnarray*}
p(3.38k,7,k,\ell)=3.764k^2-3.76k\ell+2.42k+3\ell+3>0
\end{eqnarray*}
and
\begin{eqnarray*}
q(3.38k,k,k,\ell)=1.8788k^3-1.26k^2\ell+9.9268k^2-1.74k\ell-43.72k+3\ell-15>0,
\end{eqnarray*}
\begin{eqnarray*}
q(3.38k,7,k,\ell)=15.8184k^2-13.56k\ell+16.1k+24\ell-78>0.
\end{eqnarray*}
Hence $T(n,s,k,\ell,3)>0$ for $k \geq \ell$ and $s \geq 7$.

Suppose $s=6$. Then
\begin{eqnarray*}
T(n,6,k,\ell,3)&=&\frac{(6(n-5)-3(k-3))(6(n-5)-5(\ell-5))(n-k-3)(n-\ell-1)}
{(3(n-2k)+6(k-2))(5(n-2\ell+4)+6(\ell-4))(n-5)^2}
\end{eqnarray*}
For $\ell \geq k$, noting that
\begin{eqnarray*}
&&(3(n-2k)+6(k-2)+0.619k)(5(n-2\ell+4)+6(\ell-4)-0.733k)\\&&-(3(n-2k)+6(k-2))(5(n-2\ell+4)+6(\ell-4))\\
&\geq&0.55248k\ell-0.453727k^2+6.32k>0,
\end{eqnarray*}
we have
\begin{eqnarray*}
T(n,6,k,\ell,3) \geq \frac{(6(n-5)-3(k-3)+\ell)(6(n-5)-5(\ell-5)-\ell)(n-k-3)(n-\ell-1)}
{(3(n-2k)+6(k-2)+0.619k)(5(n-2\ell+4)+6(\ell-4)-0.733k)(n-5)^2}.
\end{eqnarray*}
Since $n\geq 3.38\ell$ and $\ell \geq k$, we have
\begin{eqnarray*}
&&(n-k-3)(6(n-5)-5(\ell-5))-(n-5)(3(n-2k)+6(k-2)+0.619k)\\
&\geq&(3.38\ell-k-3)(6(3.38\ell-5)-5(\ell-5))-(3.38\ell-5)(3(3.38\ell-2k)+6(k-2)+0.619k)\\
&\geq&17.3732\ell^2+\ell(-17.37222k+28.52)+8.095k-45\\
&\geq&0.00098k^2+36.615k-45>0
\end{eqnarray*}
and
\begin{eqnarray*}
&&(n-\ell-1)(6(n-5)-3(k-3))-(n-5)(5(n-2\ell+4)+6(\ell-4)-0.733k))\\
&\geq&(2.38\ell-1)(6(3.38\ell-5)-3(k-3))-(3.38\ell-5)(5(1.38\ell+4)+6(\ell-4)-0.733k)\\
&\geq&4.6644\ell^2+\ell(-4.66246k+7.76)-0.665k+1\geq0.00194k^2+7.095k+1>0.
\end{eqnarray*}
Hence $T(n,6,k,\ell,3)>0$ for $\ell \geq k$.

For $k>\ell$, noting that
\begin{eqnarray*}
&&(3(n-2k)+6(k-2)+0.619\ell)(5(n-2\ell+4)+6(\ell-4)-0.733\ell)\\&&-(3(n-2k)+6(k-2))(5(n-2\ell+4)+6(\ell-4))\\
&=&0.896n\ell-2.929727\ell^2+6.32\ell>0,
\end{eqnarray*}
we have
\begin{eqnarray*}
T(n,6,k,\ell,3) \geq \frac{(6(n-5)-3(k-3)+\ell)(6(n-5)-5(\ell-5)-\ell)(n-k-3)(n-\ell-1)}
{(3(n-2k)+6(k-2)+0.619\ell)(5(n-2\ell+4)+6(\ell-4)-0.733\ell)(n-5)^2}
\end{eqnarray*}

Since $n\geq 3.38k$, $k\geq \ell$ and $k\geq s \geq 6$, it follows that
\begin{eqnarray*}
&&(n-k-3)(6(n-5)-5(\ell-5))-(n-5)(3(n-2k)+6(k-2)+0.619\ell)\\
&\geq&(3.38k-k-3)(6(3.38k-5)-5(\ell-5))-(3.38k-5)(3(3.38k-2k)+6(k-2)+0.619\ell)\\
&\geq&13.9932k^2+k(-13.99222\ell+18.52)+18.095\ell-45\\
&\geq&0.00098\ell^2+36.615\ell-45>0
\end{eqnarray*}
and
\begin{eqnarray*}
&&(n-\ell-1)(6(n-5)-3(k-3))-(n-5)(5(n-2\ell+4)+6(\ell-4)-0.733\ell))\\
&\geq&(3.38k-\ell-1)(6(3.38k-5)-3(k-3))-(3.38k-5)(5(3.38k-2\ell+4)+6(\ell-4)-0.733\ell)\\
&\geq&1.2844k^2+k(-1.28246\ell+9.76)-2.665\ell+1\\
&\geq&0.00194\ell^2+7.095\ell+1>0.
\end{eqnarray*}
Hence $T(n,6,k,\ell,3)>0$ for $k \geq \ell$.
Therefore, $T(n,s,k,\ell,i)>1$ for $i=3$ and $s \geq 6$.

\section{Proof of Theorem \ref{mainthm}}\label{sec:main}
\textbf{Proof.} Let $n$, $k$ and $\ell$ be  positive integers with {$n\geq 3.38\max\{k,\ell\}$}. Suppose that $\mathcal A\subseteq\binom{[n]}{k}$ and $\mathcal B\subseteq\binom{[n]}{\ell}$ are  cross-$2$-intersecting families such that $|\mathcal A||\mathcal B|$ is maximum. Then,
  \begin{eqnarray}\label{gkl}|\mathcal A||\mathcal B|\geq \binom{n-2}{k-2}\binom{n-2}{\ell-2}.\end{eqnarray} Let $g(\mathcal A)\in G_*(\mathcal A)$ and $g(\mathcal B)\in G_*(\mathcal B)$. Set $s=s^+(g(\mathcal A))$.
  Let $i$ be the minimal integer such that $g^*_i(\mathcal A)\neq \emptyset$, then $g^*_{s+2-i}(\mathcal B)\neq\emptyset$.
  To complete the proof, it suffices to verify that $g(\mathcal A)=g(\mathcal B)=\{[2]\}$. We will prove this fact through the following six scenarios.
\begin{enumerate}
\item  If $s=2$, then result holds.

 \item If $s=3$, then either $g(\mathcal A)=[3]$ and $g(\mathcal B)=\binom{[3]}{2}$ or  $g(\mathcal A)=\binom{[3]}{2}$ and $g(\mathcal B)=[3]$.
 By symmetry, we soppose that the former holds.
 Then  $|\mathcal A|=\binom{n-3}{k-3}$ and $|\mathcal B|=3\binom{n-3}{\ell-2}+\binom{n-3}{\ell-3}$.  Thus
 \begin{eqnarray*}
 |\mathcal A||\mathcal B|-\binom{n-2}{k-2}\binom{n-2}{\ell-2}&=&3\binom{n-3}{k-3}\binom{n-3}{\ell-2}+\binom{n-3}{k-3}\binom{n-3}{\ell-3}-\binom{n-2}{\ell-2}\binom{n-2}{k-2}\\
 &<& (3\binom{n-3}{k-3}-\binom{n-2}{k-2})\binom{n-2}{\ell-2}+\binom{n-3}{k-3}\binom{n-3}{\ell-3}<0
 \end{eqnarray*}
   as $\binom{n-2}{k-2}>3.38\binom{n-3}{k-3}$ and  $\binom{n-2}{\ell-2}>3.38\binom{n-3}{\ell-3}$ under the condition $n\geq 3.38\max\{k,\ell\}$, contradicting with (\ref{gkl}).
    So we may assume that $s\geq 4$.

 \item If $s=4$, noting $i\geq 3$ and $4+2-i\geq 3$, we have $i=3$. In other words, we will prove the case $(s,i)=(4,3)$. By symmetry, we may assume that $124\in g^*_3(\mathcal A)$ and $134\in g^*_3(\mathcal B)$. Then $12\not\in g(\mathcal A)$ and $g_2(\mathcal A)=\emptyset$.
 Suppose that $12\in g(\mathcal B)$. Then $g(\mathcal A)=\{123,124\}$ and the  maximality of $\mathcal B$ implies $g(\mathcal B)=\{12,134,234\}$.
Hence, $|\mathcal A|=\binom{n-3}{k-3}+\binom{n-4}{k-3}=\binom{n-2}{k-2}-\binom{n-4}{k-2}$ and $|\mathcal B|=\binom{n-2}{\ell-2}+2\binom{n-4}{\ell-3}$. Noting that  $\binom{n-4}{k-2}>2\binom{n-4}{k-3}>4\binom{n-4}{k-4}$ and $\binom{n-4}{\ell-2}>2\binom{n-4}{\ell-3}>4\binom{n-4}{\ell-4}$ under the assumption that $n\geq 3.38\max\{k,\ell\}$.
 Then
\begin{eqnarray*}
|\mathcal A||\mathcal B|&-&\binom{n-2}{k-2}\binom{n-2}{\ell-2}=2\binom{n-4}{\ell-3}\binom{n-2}{k-2}-\binom{n-4}{k-2}\binom{n-2}{\ell-2}-2\binom{n-4}{k-2}\binom{n-4}{\ell-3}\\&=&
2\binom{n-4}{\ell-3}(2\binom{n-4}{k-3}+\binom{n-4}{k-4})-\binom{n-4}{k-2}(\binom{n-4}{\ell-2}+2\binom{n-4}{\ell-3}+\binom{n-4}{\ell-4})\\&<&
5\binom{n-4}{\ell-3}\binom{n-4}{k-3}-8\binom{n-4}{k-3}\binom{n-4}{\ell-3}<0,
\end{eqnarray*}
contradicting with (\ref{gkl}).
So $12\not\in g(\mathcal B)$. Then $g_2(\mathcal A)=g_2(\mathcal B)=\emptyset$. By the maximality of $\mathcal A$ and $\mathcal B$, we get that $g(\mathcal A)=g(\mathcal B)=\binom{[4]}{3}$. In this case, $\mathcal A=\{A\in\binom{[n]}{k}: |A\cap [4]|\geq 3\}$ is a $2$-intersecting family of $\binom{[n]}{k}$, and so $|\mathcal A|<\binom{n-2}{k-2}$ holds by Theorem \ref{EKR1961} since $n\geq 3.5k$. Similarly, we can also deduce  $|\mathcal B|<\binom{n-2}{\ell-2}$. Consequently, $|\mathcal A||\mathcal B|<\binom{n-2}{k-2}\binom{n-2}{\ell-2}$, yielding a contradiction.

\item If $s=5$, then $i=3$ or $4$ since $3\leq i\leq s-1$. We prove the case $(s,i)=(5,3)$. So $g^*_3(\mathcal A)\neq\emptyset$ and hence $g^*_4(\mathcal B)\neq\emptyset$.

\textbf{Case 1}: $g_2(\mathcal{A}) = \emptyset$.

Then $125 \in g_3^*(\mathcal{A})$ and $\{1345, 2345\} \subseteq g_4^*(\mathcal{B})$. Furthermore, we can deduce that either $g_4^*(\mathcal{A}) = \{1345, 2345\}$ or $\emptyset$. If the former holds, then $g(\mathcal{A}) = g(\mathcal{B}) = \{123, 124, 125, 1345, 2345\}$. Consequently, $\mathcal{A}$ and $\mathcal{B}$ are nontrivial $2$-intersecting families of $\binom{[n]}{k}$ and $\binom{[n]}{\ell}$, respectively. Hence, $|\mathcal{A}| < \binom{n-2}{k-2}$ and $|\mathcal{B}| < \binom{n-2}{\ell-2}$, and thus $|\mathcal{A}||\mathcal{B}| < \binom{n-2}{k-2}\binom{n-2}{\ell-2}$, which is a contradiction. Now suppose the latter holds, that is, $g_4^*(\mathcal{A}) = \emptyset$. If $134 \notin g(\mathcal{A})$, then $g(\mathcal{A}) = \{123, 124, 125\}$ and $g(\mathcal{B}) = \{12, 1345, 2345\}$.

By the assumption that $n \geq 3.38k$ and direct calculations, we can obtain the following inequalities:
\[
\frac{\binom{n-3}{k-3}}{\binom{n-2}{k-2}} = \frac{k-2}{n-2} < \frac{1}{3.38} < 0.3, \quad \frac{\binom{n-4}{k-3}}{\binom{n-2}{k-2}} = \frac{(k-2)(n-k)}{(n-2)(n-3)} < 0.21,
\]
\[
\frac{\binom{n-4}{\ell-4}}{\binom{n-2}{\ell-2}} < 0.09, \quad \text{and} \quad \frac{\binom{n-5}{\ell-4}}{\binom{n-2}{k-2}} < 0.063.
\]
These elementary inequalities will be frequently used to scale the inequalities in the following proofs without further explanation.

Hence,
\begin{eqnarray*}
|\mathcal{A}||\mathcal{B}| &=& \left(\binom{n-3}{k-3} + \binom{n-4}{k-3} + \binom{n-5}{k-3}\right)\left(\binom{n-2}{\ell-2} + 2\binom{n-5}{\ell-4}\right) \\
&<& 0.72 \times 1.18 \binom{n-2}{k-2}\binom{n-2}{\ell-2} < \binom{n-2}{k-2}\binom{n-2}{\ell-2},
\end{eqnarray*}
which is a contradiction.

Now, suppose $134 \in g_3(\mathcal{A})$. It is easy to see that $234 \in g_3(\mathcal{A})$. Thus, $\{123, 124, 125, 134, 234\} \subseteq g_3(\mathcal{A})$, which implies $g_2(\mathcal{B}) = \emptyset$ and $g_3(\mathcal{B}) \subseteq \{123, 124\}$. Since $g_2(\mathcal{B}) = \emptyset$, exactly one of $123$ and $1235$ is in $g(\mathcal{B})$. A similar result holds for the pairs $124$ and $1245$.

\textbf{subase 1.1}: $123 \notin g_3(\mathcal{B})$. Then $g_3(\mathcal{B}) = \emptyset$, and hence $g(\mathcal{B}) = \binom{[5]}{4}$ and $g(\mathcal{A}) = \binom{[5]}{3}$. Thus
\begin{eqnarray*}
|\mathcal{A}||\mathcal{B}| &=& \left(\binom{n-3}{k-3} + 3\binom{n-4}{k-3} + 6\binom{n-5}{k-3}\right)\left(\binom{n-4}{\ell-4} + 4\binom{n-5}{\ell-4}\right) \\
&<& 2.19 \times 0.45 \binom{n-2}{k-2}\binom{n-2}{\ell-2} < \binom{n-2}{k-2}\binom{n-2}{\ell-2},
\end{eqnarray*}
which is a contradiction. 

\textbf{subcase 1.2}: $123 \in g_3(\mathcal{B})$. Then {$1235 \notin g_4(\mathcal{B})$}.
If  {$124 \notin g_3(\mathcal{B})$}, then $g(\mathcal{B}) = \{123, 1245, 1345, 2345\}$, and $g(\mathcal{A}) = \{123, 124, 125, 134, 234, 135, 235\}$. Thus,
\begin{eqnarray*}
|\mathcal{A}||\mathcal{B}| &=& \left(\binom{n-3}{k-3} + 3\binom{n-4}{k-3} + 3\binom{n-5}{k-3}\right)\left(\binom{n-3}{\ell-3} + 3\binom{n-5}{\ell-4}\right) \\
&<& 1.56 \times 0.57 \binom{n-2}{k-2}\binom{n-2}{\ell-2} < \binom{n-2}{k-2}\binom{n-2}{\ell-2},
\end{eqnarray*}
which is a contradiction. If  {$124 \in g_3(\mathcal{B})$}, then $g(\mathcal{B}) = \{123, 124, 1345, 2345\}$ and $g(\mathcal{A}) = \{123, 124, 125, 134, 234\}$. Thus,
\begin{eqnarray*}
|\mathcal{A}||\mathcal{B}| &=& \left(\binom{n-3}{k-3} + 3\binom{n-4}{k-3} + \binom{n-5}{k-3}\right)\left(\binom{n-3}{\ell-3} + \binom{n-4}{\ell-3} + 2\binom{n-5}{\ell-4}\right) \\
&<& 1.14 \times 0.69 \binom{n-2}{k-2}\binom{n-2}{\ell-2} < \binom{n-2}{k-2}\binom{n-2}{\ell-2},
\end{eqnarray*}
which is a contradiction.

 \textbf{Case 2}: $g_2(\mathcal A)=\{12\}$.

 Then $[2]\subset E$ for all $E\in g(\mathcal B)$. By the assumption that $g_3^*(\mathcal A)\neq\emptyset$ and $13\not\in g_2(\mathcal A)$, it follows that  $\{135,235\}\subseteq g_3^*(\mathcal A)$ and $1245\in g^*_4(\mathcal B)$. Then either $123 \in g_3(\mathcal B)$ or $1235 \in g^*_4(\mathcal B)$. If the former holds,   it is easy to get that $g(\mathcal A)=\{12, 134, 135, 234, 235\}$ and $g(\mathcal B)=\{123, 1245\}$,  and then
    \begin{eqnarray*}
    |\mathcal A||\mathcal B|&=&(\binom{n-2}{k-2}+2\binom{n-4}{k-3}+2\binom{n-5}{k-3})(\binom{n-3}{\ell-3}+\binom{n-5}{\ell-4})\\&<&1.84\times0.39\binom{n-2}{k-2}\binom{n-2}{\ell-2}<\binom{n-2}{k-2}\binom{n-2}{\ell-2},
\end{eqnarray*}
a contradiction.
   If the latter holds, then $g(\mathcal A)=\{12, 134, 135, 234, 235,  {145}, 245, 345\}$ and $g(\mathcal B)=\{1234, 1235, 1245\}$, and so we can deduce that
      \begin{eqnarray*}
    |\mathcal A||\mathcal B|&=&(\binom{n-2}{k-2}+2\binom{n-4}{k-3}+ {5}\binom{n-5}{k-3})(\binom{n-4}{\ell-4}+2\binom{n-5}{\ell-4})\\
    &<& {2.47}\times0.27\binom{n-2}{k-2}\binom{n-2}{\ell-2}<\binom{n-2}{k-2}\binom{n-2}{\ell-2},
\end{eqnarray*}
a contradiction.

\textbf{Case 3}: $g_2(\mathcal A)=\{12,13,23\}$.

 It is easy to verify that $g(\mathcal A)=\{12,13,23,145,245,345\}$ and $g(\mathcal B)=\{1234,1235\}$. Then,
\begin{eqnarray*}
|\mathcal A||\mathcal B|&=&(\binom{n-2}{k-2}+2\binom{n-3}{k-2}+3\binom{n-5}{k-3})(\binom{n-4}{\ell-4}+\binom{n-5}{\ell-4})\\&<&4\times 0.18\binom{n-2}{k-2}\binom{n-2}{\ell-2}<\binom{n-2}{k-2}\binom{n-2}{\ell-2},
\end{eqnarray*}
a contradiction.

\item In the following, we assume $s\geq 5$ and $(s,i)\neq (6,4)$.

 Set $\mathcal A_1=\mathcal A\cup\mathcal D(g_i^*(\mathcal A)')$ and  $\mathcal B_1=\mathcal B\backslash\mathcal D(g^*_{s+2-i}(\mathcal B))$. By (iii) of Lemma \ref{lemmag1}, $\mathcal A_1$ and $\mathcal B_1$ are also cross-$2$-intersecting with
$|\mathcal A_1|=|\mathcal A|+|g^*_i(\mathcal A)|\binom{n-s}{k-i+1}$ and $|\mathcal B_1|=|\mathcal B|-|g^*_{s+2-i}(\mathcal B)|\binom{n-s}{\ell-(s+2-i)}$.
Since $|\mathcal A||\mathcal B|$ is maximum, we have $|\mathcal A_1||\mathcal B_1|\leq |\mathcal A||\mathcal B|$, it follows that
\begin{equation}
\label{ineq31}
\frac{|\mathcal B|}{|\mathcal A|+|g^{*}_{i}(\mathcal A)|\binom {n-s}{k-i+1}}
\leq \frac{|g^{*}_{s+2-i}(\mathcal B)|\binom {n-s}{\ell-s-2+i}}{|g^{*}_{i}(\mathcal A)|\binom {n-s}{k-i+1}}.
\end{equation}
 Set $\mathcal A_2=\mathcal A\backslash\mathcal D(g_i^*(\mathcal A))$ and  $\mathcal B_2=\mathcal B\cup\mathcal D(g^*_{s+2-i}(\mathcal B)')$. Analogously, $\mathcal A_1$ and $\mathcal B_1$ are also cross-$2$-intersecting with
$|\mathcal A_2|=|\mathcal A|-|g^*_i(\mathcal A)|\binom{n-s}{k-i+1}$ and $|\mathcal B_2|=|\mathcal B|-|g^*_{s+2-i}(\mathcal B)|\binom{n-s}{\ell-(s+2-i)}$.
Moreover, we can get in the similar way that
\begin{equation}
\label{ineq32}
\frac{|\mathcal A|}{|\mathcal B|+|g^{*}_{s+2-i}(\mathcal B)|\binom {n-s}{\ell-(s+2-i)+1}}
\leq \frac{|g^{*}_{ {i}}(\mathcal A)|\binom {n-s}{k-i}}{|g^{*}_{s+2-i}(\mathcal B)|\binom {n-s}{\ell-(s+2-i)+1}}.
\end{equation}
Combining with  (\ref{ineq31}) and (\ref{ineq32}), we obtain
\begin{eqnarray}
\label{ineq33}
\frac{|\mathcal A|} {|\mathcal A|+|g^*_i(\mathcal A)|\binom {n-s}{k-i+1}}
\frac{|\mathcal B|}{|\mathcal B|+|g^*_{s+2-i}(\mathcal B)|\binom {n-s}{\ell-s-2+i}}
\leq \frac{\binom {n-s}{k-i} \binom {n-s}{\ell-s-2+i}}{\binom {n-s}{k-i+1}\binom {n-s}{\ell-s-1+i}}.
\end{eqnarray}
Set $\nabla(g_i^*(\mathcal A)')=\{F\in\binom{[s-1]}{i}: \mbox{$E\subset F$ for some $E\in g_i^*(\mathcal A) {'}$}\}$. By Sperner Theorem in \cite{Sperner},  we have
$$\frac{|\nabla(g_i^*(\mathcal A)')|}{|g_i^*(\mathcal A)|}\geq \frac{\binom{s-1}i}{\binom{s-1}{i-1}}. $$
Moreover, $\mathcal D(\nabla(g_i^*(\mathcal A)'))\subset \mathcal A$. Therefore,\[|\mathcal A|\geq |\mathcal D(\nabla(g_i^*(\mathcal A)')\cup g_i^*(\mathcal A))|\geq   \frac{|g^*_i(\mathcal A)|}{\binom{s-1}{i-1}}(\binom{s-1}{i-1}\binom{n-s}{k-i}+\binom{s-1}{i}\binom{n-s+1}{k-i})\]
Note that the ratio $\frac{|\mathcal A|} {|\mathcal A|+|g^*_i(\mathcal A)|\binom {n-s}{k-i+1}}$ increase with $|\mathcal A|$. It follows that
\begin{equation}\label{ab1}\frac{|\mathcal A|}{|\mathcal A|+|g^*_i(\mathcal A)|\binom{n-s}{k-i+1}}\geq\frac{\binom{s-1}{i-1}\binom{n-s}{k-i}+\binom{s-1}{i}\binom{n-s+1}{k-i}}
{\binom{s-1}{i-1}\binom{n-s+1}{k-i+1}+\binom{s-1}{i}\binom{n-s+1}{k-i}}
= {f(n,k,s,i)}.
 \end{equation}
 In the similar way, we can also deduce
 \begin{eqnarray}\nonumber\frac{|\mathcal B|}{|\mathcal B|+|g^*_{s+2-i}(\mathcal B)|\binom{n-s}{\ell-(s+2 -i)+1}}&\geq&\frac{\binom{s-1}{s+1-i}\binom{n-s}{\ell-(s+2-i)}+\binom{s-1}{s+2-i}\binom{n-s+1}{\ell-(s+2-i)}}
{\binom{s-1}{s+1-i}\binom{n-s+1}{\ell-(s+2-i)+1}+\binom{s-1}{s+2-i}\binom{n-s+1}{\ell-(s+2-i)}}\\ \label{ab2}&=& {f(n,\ell,s,s+2-i)}.
 \end{eqnarray}
 Combining with (\ref{ineq33}), (\ref{ab1}) and (\ref{ab2}),
 we have
\begin{eqnarray*} {f(n,k,s,i)f(n,\ell,s,s+2-i)}&\leq& \frac{|\mathcal A|} {|\mathcal A|+|g^*_i(\mathcal A)|\binom {n-s}{k-i+1}}
\frac{|\mathcal B|}{|\mathcal B|+|g^*_{s+2-i}(\mathcal B)|\binom {n-s}{\ell-s-2+i}}\\ &\leq& \frac{\binom {n-s}{k-i} \binom {n-s}{\ell-s-2+i}}{\binom {n-s}{k-i+1}\binom {n-s}{\ell-s-1+i}}\\&=&\frac{(k-i+1)(\ell+i-s-1)}{(n+i-k-s)(n+2-\ell-i)},\end{eqnarray*}
  contracting with Lemma \ref{lem:main}.
  \item $(s,i)=(6,4)$. Let
$g(\mathcal A)\in G_*(\mathcal A)$ and $g(\mathcal B)\in G_*(\mathcal B)$. Suppose $g^*_4(\mathcal A)\neq\emptyset$ and $g^*_4(\mathcal B)\neq\emptyset$. We distinguish two cases to complete the proof. {In order to make this section readable, we will only prove $i=0,1$ in case 1 and only list all the specific cases in other scenarios, since the proofs are all similar with the case  $(s,i)=(5,3)$ which will be omitted.}

{\bf Case 4: $g_2(\mathcal A)\cup g_2(\mathcal B)=\emptyset$. }

Set $|g_3(\mathcal A)\cup g_3(\mathcal B)|=i$ and suppose  $|g_3(\mathcal A)|\geq |g_3(\mathcal B)|$.

 If $i=0$, then $g^*_4(\mathcal A)=g^*_4(\mathcal B)=\binom{[6]}{4}$, and  $\mathcal A=\mathcal F(n,k,2,2)$ and $\mathcal B=\mathcal F(n,\ell,2,2)$. Hence $|\mathcal A||\mathcal B|<\binom{n-2}{k-2}\binom{n-2}{\ell-2}$.

 If $i=1$, then $g_3(\mathcal A)=\{123\}$ and $g_3(\mathcal B)\subseteq \{123\}$. If $g_3(\mathcal B)=\{123\}$, we have that  $$g(\mathcal A)=g(\mathcal B)=\{123,1245,1246,1256,1345,2345,1346,2346,1356,2356\},$$ then $\mathcal A$ and $\mathcal B$ are two nontrivial $2$-intersecting families of $\binom{[n]}{k}$ and $\binom{[n]}{\ell}$, respectively. Hence $|\mathcal A||\mathcal B|<\binom{n-2}{k-2}\binom{n-2}{\ell-2}$. If  $g_3(\mathcal B)=\emptyset$, we have
 $$g(\mathcal A)=\{123,1245,1246, {1256},1345,2345,1346,2346,1356,2356,1456,2456,3456\}$$ and $$g(\mathcal B)=\{ {1234}, {1235}, {1236},1245,1246, {1256},1345,2345,1346,2346,1356,2356\}.$$
Consequently,
\begin{eqnarray*}
 {|\mathcal A||\mathcal B|}&=&(\binom{n-3}{k-3}+3\binom{n-5}{k-4}+9\binom{n-6}{k-4})(\binom{n-4}{\ell-4}+4\binom{n-5}{\ell-4}+7\binom{n-6}{\ell-4})\\
&<&(0.3+12\times 0.063)\times (0.09+11\times 0.063)\binom{n-2}{k-2})\binom{n-2}{\ell-2}<\binom{n-2}{k-2}\binom{n-2}{\ell-2}.
\end{eqnarray*}
 If $i=2$, then $g_3(\mathcal A)=\{123,124\}$ and $g_3(\mathcal B)\subseteq \{123,124\}$.  Hence either $$g(\mathcal A)=g(\mathcal B)=\{123,124,1256,1345,2345,1346,2346\},$$ or
 $$g(\mathcal A)=\{123,124,1256,1345,2345,1346,2346,1356,2356\}$$ and $$g(\mathcal B)=\{123,1245,1246, {1256},1345,2345,1346,2346\},$$

 Then, \begin{eqnarray*}
 {|\mathcal A||\mathcal B|}&=&(\binom{n-3}{k-3}+\binom{n-4}{k-3}+2\binom{n-5}{k-4}+5\binom{n-6}{k-4})(\binom{n-3}{\ell-3}+3\binom{n-5}{\ell-4}+3\binom{n-6}{\ell-4})\\
&<&(0.3+0.21+7\times 0.063)\times (0.3+6\times 0.063)\binom{n-2}{k-2})\binom{n-2}{\ell-2}<\binom{n-2}{k-2}\binom{n-2}{\ell-2}.
\end{eqnarray*}

or $$g(\mathcal A)=\{123,124,1256,1345,2345,1346,2346,2356,1456,2456,3456\}$$
and $$g(\mathcal B)=\{1234,1235,1236,1245,1246, {1256},1345,2345,1346,2346\}.$$

 Then, \begin{eqnarray*}
 {|\mathcal A||\mathcal B|}&=&(\binom{n-3}{k-3}+\binom{n-4}{k-3}+2\binom{n-5}{k-4}+7\binom{n-6}{k-4})(\binom{n-4}{\ell-4}+4\binom{n-5}{\ell-4}+5\binom{n-6}{\ell-4})\\
&<&(0.3+0.21+9\times 0.063)\times (0.09+9\times 0.063)\binom{n-2}{k-2})\binom{n-2}{\ell-2}<\binom{n-2}{k-2}\binom{n-2}{\ell-2}.
\end{eqnarray*}

If $i=3$, then  $g_3(\mathcal A)=\{123,124,125\}$ and  $g_3(\mathcal B)\subseteq\{123,124,125\}$. Therefore,
 either  $$g(\mathcal A)=\{123,124,125,1345,2345,1346,2346\} \:\: \text{and}\:\: g(\mathcal B)=\{123,124,1256,1345,2345\},$$
  Then, \begin{eqnarray*}
 {|\mathcal A||\mathcal B|}&=&(\binom{n-3}{k-3}+\binom{n-4}{k-3}+\binom{n-5}{k-3}+2\binom{n-5}{k-4}+2\binom{n-6}{k-4})(\binom{n-3}{\ell-3}+\binom{n-4}{\ell-3}+2\binom{n-5}{\ell-4}+\binom{n-6}{\ell-4})\\
&<&(0.3+2\times0.21+4\times 0.063)\times (0.3+0.21+3\times 0.063)\binom{n-2}{k-2})\binom{n-2}{\ell-2}<\binom{n-2}{k-2}\binom{n-2}{\ell-2}.
\end{eqnarray*}

 or   $$g(\mathcal A)=\{123,124,125,1345,2345,1346,2346,1356,2356\}$$ and $$g(\mathcal B)=\{123,1245,1246,1256,1345,2345\},$$
   Then, \begin{eqnarray*}
 {|\mathcal A||\mathcal B|}&=&(\binom{n-3}{k-3}+\binom{n-4}{k-3}+\binom{n-5}{k-3}+2\binom{n-5}{k-4}+4\binom{n-6}{k-4})(\binom{n-3}{\ell-3}+3\binom{n-5}{\ell-4}+2\binom{n-6}{\ell-4})\\
&<&(0.3+2\times0.21+6\times 0.063)\times (0.3+5\times 0.063)\binom{n-2}{k-2})\binom{n-2}{\ell-2}<\binom{n-2}{k-2}\binom{n-2}{\ell-2}.
\end{eqnarray*}

or $$g(\mathcal A)=\{123,124,125,1345,2345,1346,1356,2346,2356,1456,2456,3456\}$$
and $$g(\mathcal B)=\{1234,1235,1236,1245,1246,1256,1345,2345\}.$$
  Then, \begin{eqnarray*}
 {|\mathcal A||\mathcal B|}&=&(\binom{n-3}{k-3}+\binom{n-4}{k-3}+\binom{n-5}{k-3}+2\binom{n-5}{k-4}+7\binom{n-6}{k-4})(\binom{n-4}{\ell-4}+4\binom{n-5}{\ell-4}+3\binom{n-6}{\ell-4})\\
&<&(0.3+2\times0.21+9\times 0.063)\times (0.09+7\times 0.063)\binom{n-2}{k-2})\binom{n-2}{\ell-2}<\binom{n-2}{k-2}\binom{n-2}{\ell-2}.
\end{eqnarray*}

If $i=4$, then $g_3(\mathcal A)=\{123,124,134,234\}$ and  $g_3(\mathcal B)\subseteq\{123,124,134,234\}$. Therefore, either
 $$g(\mathcal A)=\{123,124,134,234,1256,1356,1456,2356,2456,3456\}$$ and
 $$g(\mathcal B)=\{1236, {1235},1246, {1245},1346, {1345},2346, {2345}\},$$
   Then, \begin{eqnarray*}
 {|\mathcal A||\mathcal B|}&=&(\binom{n-3}{k-3}+3\binom{n-4}{k-3}+5\binom{n-6}{k-4})(4\binom{n-5}{\ell-4}+4\binom{n-6}{\ell-4})\\
&<&(0.3+3\times0.21+5\times 0.063)\times8\times 0.063\binom{n-2}{k-2})\binom{n-2}{\ell-2}<\binom{n-2}{k-2}\binom{n-2}{\ell-2}.
\end{eqnarray*}

or $$g(\mathcal A)=\{123,124,134,234,1256,1356,2356\} \:\: \text{and}\:\: g(\mathcal B)=\{123, {1245},1246, {1345},1346, {2345},2346\},$$
   Then, \begin{eqnarray*}
 {|\mathcal A||\mathcal B|}&=&(\binom{n-3}{k-3}+3\binom{n-4}{k-3}+3\binom{n-6}{k-4})(\binom{n-3}{\ell-3}+3\binom{n-5}{\ell-4}+3\binom{n-6}{\ell-4})\\
&<&(0.3+3\times0.21+3\times 0.063)\times(0.3+6\times 0.063)\binom{n-2}{k-2})\binom{n-2}{\ell-2}<\binom{n-2}{k-2}\binom{n-2}{\ell-2}.
\end{eqnarray*}

or $$g(\mathcal A)=\{123,124,134,234,1256\}\:\: \text{and}\:\: g(\mathcal B)=\{123,124, {1345},1346, {2345},2346\}.$$
   Then, \begin{eqnarray*}
 {|\mathcal A||\mathcal B|}&=&(\binom{n-3}{k-3}+3\binom{n-4}{k-3}+\binom{n-6}{k-4})(\binom{n-3}{\ell-3}+\binom{n-4}{\ell-3}+2\binom{n-5}{\ell-4}+2\binom{n-6}{\ell-4})\\
&<&(0.3+3\times0.21+0.063)\times(0.3+0.21+4\times 0.063)\binom{n-2}{k-2})\binom{n-2}{\ell-2}<\binom{n-2}{k-2}\binom{n-2}{\ell-2}.
\end{eqnarray*}

If $|g_3(\mathcal A)|=5$, then $g_3(\mathcal A)=\{123,124,134,234,125\}$, in this case, it is easy to find that $g(\mathcal B)\subseteq \{1235,1245,1345,2345\}$, contradicting that $g_4^*(\mathcal B)\neq\emptyset$.

\smallskip
{\bf Case 5: $g_2(\mathcal A)\cup g_2(\mathcal B)\neq\emptyset$. }

Suppose $g_2(\mathcal A)\neq \emptyset$. If  $g_2(\mathcal B)\neq \emptyset$, it is easily seen that  $ g(\mathcal A)=g(\mathcal B)=\{12\}$, contradicting that $g^*_4(\mathcal A)=\emptyset$. Therefore, $g_2(\mathcal B)=\emptyset$.
 We distinguish some cases to complete this proof.

\textbf{Subcase 5.1}: $g_2(\mathcal A)=\{12\}$.  If $g_3(\mathcal A)\neq\emptyset$, then $\{134,234\}\subseteq g_3(\mathcal A)$, and so $[2]\subset E$ and $|E\cap [4]|\geq 3$ for all $E\in g(\mathcal B)$. Therefore,
$$g(\mathcal A)=\{12,134,234,1356,2356,1456,2456,3456\}\:\: \text{and}\:\: g(\mathcal B)=\{1234, 1235, 1236, 1245,1246\},$$
   Then, \begin{eqnarray*}
 {|\mathcal A||\mathcal B|}&=&(\binom{n-2}{k-2}+2\binom{n-4}{k-3}+5\binom{n-6}{k-4})(\binom{n-4}{\ell-4}+2\binom{n-5}{\ell-4}+2\binom{n-6}{\ell-4})\\
&<&(1+2\times0.21+5\times0.063)\times(0.09+4\times 0.063)\binom{n-2}{k-2})\binom{n-2}{\ell-2}<\binom{n-2}{k-2}\binom{n-2}{\ell-2}.
\end{eqnarray*}

or
$$g(\mathcal A)=\{12,134,234,1356,2356\}\:\: \text{and}\:\: g(\mathcal B)=\{123,  1245,1246\},$$
   Then, \begin{eqnarray*}
 {|\mathcal A||\mathcal B|}&=&(\binom{n-2}{k-2}+2\binom{n-4}{k-3}+2\binom{n-6}{k-4})(\binom{n-3}{k-3}+\binom{n-5}{\ell-4}+\binom{n-6}{\ell-4})\\
&<&(1+2\times0.21+2\times0.063)\times(0.3+2\times 0.063)\binom{n-2}{k-2})\binom{n-2}{\ell-2}<\binom{n-2}{k-2}\binom{n-2}{\ell-2}.
\end{eqnarray*}

or $$g(\mathcal A)=\{12,134,234,135,235,1456,2456,3456\}\:\: \text{and}\:\: g(\mathcal B)=\{1234, 1235, 1236, 1245\}.$$
 Then, \begin{eqnarray*}
 {|\mathcal A||\mathcal B|}&=&(\binom{n-2}{k-2}+2\binom{n-4}{k-3}+2\binom{n-5}{k-3}+3\binom{n-6}{k-4})(\binom{n-4}{\ell-4}+2\binom{n-5}{\ell-4}+\binom{n-6}{\ell-4})\\
&<&(1+4\times0.21+5\times0.063)\times(0.09+3\times 0.063)\binom{n-2}{k-2})\binom{n-2}{\ell-2}<\binom{n-2}{k-2}\binom{n-2}{\ell-2}.
\end{eqnarray*}

Suppose $g_3(\mathcal A)=\emptyset$. Then either
\[g(\mathcal A)=\{12, {1345},1346, {2345},2346,1356,2356,1456,2456,3456\}\]and \[g(\mathcal B)=\{1234, 1235, 1236, 1245,1246,1256\},\]

 Then, \begin{eqnarray*}
 {|\mathcal A||\mathcal B|}&=&(\binom{n-2}{k-2}+2\binom{n-5}{k-4}+7\binom{n-6}{k-4})(\binom{n-4}{\ell-4}+2\binom{n-5}{\ell-4}+3\binom{n-6}{\ell-4})\\
&<&(1+9\times0.063)\times(0.09+5\times 0.063)\binom{n-2}{k-2})\binom{n-2}{\ell-2}<\binom{n-2}{k-2}\binom{n-2}{\ell-2}.
\end{eqnarray*}

or $$g(\mathcal A)=\{12, {1345},1346, {2345},2346,1356,2356\}\:\: \text{and}\:\: g(\mathcal B)=\{123,1245,1246,1256\},$$
 Then, \begin{eqnarray*}
 {|\mathcal A||\mathcal B|}&=&(\binom{n-2}{k-2}+2\binom{n-5}{k-4}+4\binom{n-6}{k-4})(\binom{n-3}{\ell-3}+\binom{n-5}{\ell-4}+2\binom{n-6}{\ell-4})\\
&<&(1+6\times0.063)\times(0.3+3\times 0.063)\binom{n-2}{k-2})\binom{n-2}{\ell-2}<\binom{n-2}{k-2}\binom{n-2}{\ell-2}.
\end{eqnarray*}

or $$g(\mathcal A)=\{12, {1345},1346, {2345},2346\}\:\: \text{and}\:\: g(\mathcal B)=\{123,  124, 1256\}.$$
 Then, \begin{eqnarray*}
 {|\mathcal A||\mathcal B|}&=&(\binom{n-2}{k-2}+2\binom{n-5}{k-4}+2\binom{n-6}{k-4})(\binom{n-3}{\ell-3}+\binom{n-4}{\ell-3}+\binom{n-6}{\ell-4})\\
&<&(1+4\times0.063)\times(0.3+0.21+ 0.063)\binom{n-2}{k-2})\binom{n-2}{\ell-2}<\binom{n-2}{k-2}\binom{n-2}{\ell-2}.
\end{eqnarray*}

\textbf{Subcase 5.2}: $g_2(\mathcal A)=\binom{[3]}{2}$. Then, $[3]\subseteq E$ for all $E\in g(\mathcal B)$, and  the assumption that $g^*_4(\mathcal B)\neq\emptyset$ implies $g(\mathcal B)=\{1234,1235,1236\}$.
Therefore, $g(\mathcal A)=\{12,13,23,1456, 2456, 3456\}$.
\end{enumerate}\qed

\section{A remark on nontrivial cross-2-intersecting families}\label{sec:nontrivial}
Actually, we can obtain more detailed results when $\mathcal A\subseteq\binom{[n]}{k}$ and $\mathcal B\subseteq\binom{[n]}{\ell}$ are nontrivial cross-$2$-intersecting families, i.e.,  $|\cap_{A\in\mathcal A\cup\mathcal B}A|<2$. Before stating the result, we first define some families of sets related.
For nonnegative integer $r$, set
\[
\mathcal F_{n,k,2,r}=\left\{A\in\binom{[n]}{k}: |A\cap [2+2r]|\geq 2+r\right\}.\]

For $3\leq s\leq \min\{k,\ell\}$, put
$$\mathcal A_{n,\ell,s,2}=\{A\in\binom{[n]}{\ell}:\mbox{ $[s]\subset A$}\},$$
  $$\mathcal B_{n,k,s,2}=\{A\in\binom{[n]}{k}:\mbox{$|A\cap [s]|\geq 2$}\},$$
$$\mathcal H_{n,k,s,2}=\{A\in\binom{[n]}{k}: \mbox{either $[2]\subset A$ or $|A\cap[s]|\geq s-1$}\}.$$

 $$\mathcal I_{n,\ell,s,2}=\{A\in\binom{[n]}{\ell}:\mbox{ $[2]\subset A$ and $|A\cap [s]|\geq 3$}\}.$$
  It is obvious that $\mathcal A_{n,\ell,3,2}=\mathcal I_{n,\ell,3,2}$ and $\mathcal B_{n,k,3,2}=\mathcal H_{n,k,3,2}$.
 And we also set
  \begin{eqnarray}\label{h}
  h_{n,k,\ell,s,2}=(\binom{n-2}{k-2}+t\binom{n-s}{k-s+1})(\binom{n-2}{\ell-2}-\binom{n-s}{\ell-2})
  \end{eqnarray}
  and
  $$f_{n,k,\ell,2,r}=|\mathcal F_{n,k,2,r}||\mathcal F_{n,\ell,2,r}|.$$


The technique in~\cite{He2024} can be adapted easily, therefore the detailed proof of the following result is omitted.
\begin{theorem}\label{thm:non}Let $n, k,\ell$ be positive integers with $n\ge 3.38\max\{k,\ell\}$. If  $\mathcal A\subseteq\binom{[n]}{k}$ and $\mathcal B\subseteq\binom{[n]}{\ell}$ are nontrivial cross-$2$-intersecting families of $\binom{[n]}{k}$, i.e, $|\bigcap_{A\in\mathcal A\cup\mathcal B}A|<2$, then
 $$|\mathcal A||\mathcal B|\leq \max\left\{h_{n,k,\ell,3,2}, h_{n,k,\ell,3,t},f_{n,k,\ell,2,1}\right\}.$$
 Moreover, the equality holds if and only if, up to isomorphism, one of the following holds:
 \begin{itemize}
   \item [\rm(i)]  $\mathcal A=\mathcal H_{n,k,3,2}$ and $\mathcal B=\mathcal I_{n,\ell,3,2}$ when $h_{n,k,\ell,3,2}\geq \max\left\{h_{n,k,\ell,k+1,2},f_{n,k,\ell,2,1}\right\}$;
   \item [\rm(ii)] $\mathcal A=\mathcal H_{n,k,k+1,2}$ and $\mathcal B=\mathcal I_{n,\ell,k+1,2}$ when
    $h_{n,k,\ell,k+1,2}\geq \max\left\{ h_{n,k,\ell,3,2},f_{n,k,\ell,2,1}\right\}$;
   \item [\rm(iii)]$\mathcal A=\mathcal F_{n,k,2,1}$, $\mathcal B=\mathcal F_{n,\ell,2,1}$ when $f_{n,k,\ell,2,1}\geq\max\{ h_{n,k,\ell,3,2},h_{n,k,\ell,k+1,2}\}$.
 \end{itemize}
 \end{theorem}

\end{document}